%%%%%%%%%%%%%%%%%%%%%%%%%%%%%%%%%%%%%%%%%%%%%%%%    
%
%        THIS IS A  PLAIN TeX FILE
%
%%%%%%%%%%%%%%%%%%%%%%%%%%%%%%%%%%%%%%%%%%%%%%%%

\magnification=1200

\font\titfont=cmr10 at 12 pt

\font\headfont=cmr10 at 12 pt

%\font\AAA=Times.dfont  at 12pt
 %\font\BBB=Times.dfont at 8pt

%\font\AAA=cmr10 at 12pt
%\font\BBB=cmr10 at 8pt

\def\AAA{\bf}
\def\BBB{\bf}

\overfullrule=0in

\def\boxit#1{\hbox{\vrule
 \vtop{%
  \vbox{\hrule\kern 2pt %
     \hbox{\kern 2pt #1\kern 2pt}}%
   \kern 2pt \hrule }%
  \vrule}}

  \def\harr#1#2{\ \smash{\mathop{\hbox to .3in{\rightarrowfill}}\limits^{\scriptstyle#1}_{\scriptstyle#2}}\ }

 \def\Iff{\qquad\iff\qquad}
 \def\GG{{{\bf G} \!\!\!\! {\rm l}}\ }

\def\bra#1#2{\langle #1, #2\rangle}

\def\ss{\subset}

\def\half{\hbox{${1\over 2}$}}
\def\smfrac#1#2{\hbox{${#1\over #2}$}}

\def\dist{{\rm dist}}

\def\Sym{{\rm Sym}^2}

\def\arr{\longrightarrow}

\def\rn{\bbr^n}

\def\Symn{{\Sym(\rn)}}

\def\Theorem#1{\medskip\noindent {\bf THEOREM \bf #1.}}
\def\Prop#1{\medskip\noindent {\bf Proposition #1.}}
\def\Cor#1{\medskip\noindent {\bf Corollary #1.}}
\def\Lemma#1{\medskip\noindent {\bf Lemma #1.}}
\def\Remark#1{\medskip\noindent {\bf Remark #1.}}

\def\Def#1{\medskip\noindent {\bf Definition #1.}}

\def\pf{\medskip\noindent {\bf Proof.}\ }
\def\qed{\hfill  $\vrule width5pt height5pt depth0pt$}

   \def\cc{{\cal C}}

\def\ch{{\cal H}}

\def\vf{\varphi}

\def\wh{\widehat}

\def\and{\qquad {\rm and} \qquad}
\def\arr{\longrightarrow}

\def\bbr{{\bf R}}

\def\a{\alpha}
\def\b{\beta}
\def\d{\delta}
\def\e{\epsilon}

\def\l{\lambda}

\def\x{\xi}

\def\CH{{\cal H}}

\def\Symn{\Sym(\rn)}
 
\def\USC{{\rm USC}}
\def\fa{{\rm\ \  for\ all\ }}

\def\AA{1}
\def\BB{2}
\def\CC{3}
\def\DD{4}
\def\EE{5}
\def\FF{6}
\def\GG{7}
\def\HH{8}

\def\AAA{1}
 \def\BBB{2} 
 \def\AA{3}
\def\BB{4} 
\def\CC{5}
\def\DD{6}  
\def\EE{7} 
\def\FF{8}
\def\GG{9}
\def\HH{10}

\centerline
{
\titfont NOTES ON THE DIFFERENTIATION OF }

\smallskip

\centerline
{
\titfont QUASI-CONVEX FUNCTIONS }
\smallskip

\bigskip

\centerline{\titfont F. Reese Harvey and H. Blaine Lawson, Jr.$^*$}
\vglue .9cm
\smallbreak\footnote{}{ $ {} \sp{ *}{\rm Partially}$  supported by
the N.S.F. } 

\vskip .2in

\centerline{\bf ABSTRACT} \medskip
  \font\abstractfont=cmr10 at 10 pt

  {{\parindent= .93in
\narrower\abstractfont \noindent
This expository paper presents elementary proofs of  four  basic results concerning 
derivatives of quasi-convex functions.  They are combined into a fifth
 theorem which is simple to apply and adequate in many cases.
Along the way we establish the
 equivalence of the basic lemmas of Jensen and Slodkowski.
 
}}

\vskip.5in

%%%%%%%%%%%%%%%%%%%%%%%%%%%%%%%%%%%%%%%%%%%%%%%%%%%%
%%%%%%%%%%%%%%%%%%%%%%%%%%%%%%%%%%%%%%%%%%%%%%%%%%%%
%%%%%%%%%%%%%%%%%%%%%%%%%%%%%%%%%%%%%%%%%%%%%%%%%%%%
%%%%%%%%%%%%%%%%%%%%%%%%%%%%%%%%%%%%%%%%%%%%%%%%%%%%
%%%%%%%%%%%%%%%%%%%%%%%%%%%%%%%%%%%%%%%%%%%%%%%%%%%%

\centerline{\headfont \AAA.\ Introduction.}
 \medskip 

Let $u(x)$ be a real-valued function on an open set $X\ss\rn$.
Then $u(x)$ is said to be {\sl quasi-convex}
if the function $u(x) +{\l \over 2} |x|^2$ is convex for some $\l  \geq0$.
There are four basic results concerning the differentiability 
of such functions.  To state some of them we need the following concept.
Let $\Symn$ denote the set of $n\times n$ symmetric matrices.

\Def{\AAA.1}  
A point $x\in X$ is called an {\bf upper contact point for $u$}
if there exists $(p,A) \in\rn\times\Symn$ such that
$$
u(y)\ \leq\ u(x) + \bra p {y-x} +\half \bra{A(y-x)}{y-x} \qquad\forall\, y \ {\rm near}\  x.
\eqno{(\AAA.1)}
$$
In this case, $(p,A)$ is called an {\bf upper contact jet for $u$ at $x$}.

The first result is the differentiability at upper contact points.

\Lemma{\AAA.2. (D at UCP)} {\sl Suppose $u$ is quasi-convex.  If $x$ is an upper contact point
for $u$, then $u$ is differentiable at $x$. Moreover, if $(p,A)$ is any upper contact
jet for $u$ at $x$, then $p=D_{x}$ is unique.
}
\medskip

Another even more standard result is called partial continuity of the gradient, or first derivative.

\Lemma{\AAA.3. (PC of FD)} {\sl Suppose $u$ is quasi-convex and $x_j \to x$.
If $u$ is differentiable at each $x_j$ and at $x$, then $D_{x_j}u\to D_{x}u$.
}

\medskip

The next two results concern the second-order  contact of quasi-convex functions
and are of a deeper nature.  

\Theorem {\AAA.4. (Alexandrov)} {\sl
 A locally quasi-convex function is twice differentiable almost everywhere.
}
\medskip

For the next result we need two variations of the notion of an upper contact jet.
First, we say that $(p,A)$ is a {\bf strict} upper contact jet for $u \in \USC(X)$ at $x_0\in X$ if
the upper contact inequality (\AAA.1) is strict for $y\neq x_0$.
An understanding of the strict  upper contact jets will be  adequate for 
our discussion since $(p,A)$ is an upper contact jet 
 if and only if   $(p, A+\e I)$ is a strict upper contact jet for all $\e>0$.
Second, we need a   notion of upper contact point and jet, which requires the inequality
(\AA.1) to hold globally.

\Def{\AAA.5}  Given $u\in \USC(X)$ and $A\in \Symn$, a point $x$ is called a
{\bf global upper contact point of type $A$ on $X$} if for some $p \in\rn$
$$
u(y)\ \leq\ u(x) + \bra p {y-x}  + \half\bra{A(y-x)}{y-x}\quad \forall\, y\in X.
\eqno{(\AAA.2)}
$$
Let $\cc(u,X,A)$ denote the set of all global upper contact points of type $A$ on $X$ for
the function $u$.

\Remark{\AAA.6} Note that if $u$ is quasi-convex, then by (D at UCP) each point
$x\in \cc(u,X,A)$ is a point of differentiability and the only $p$ in (\AAA.2) is $p =D_xu$.

\Theorem{\AAA.7. (Jensen-Slodkowski)} {\sl
Suppose that $u$ is a quasi-convex function possessing a strict upper contact
jet $(p,A)$ at $x$.  Let $B_\rho$ denote the ball of radius $\rho$ about $x$.
Then there exists $\bar \rho >0$ such that the measure }
$$
\left|   \cc(u, B_\rho, A)      \right| \ >\ 0 \qquad \forall\, 0 < \rho\leq \bar \rho.
\eqno{(\AAA.3)}
$$

This result follows in a straightforward/elementary manner (see Section \BB) from
Slodkowski's Lemma \BB.1 below, which in turn is proved  in Sections \CC--\EE.

On the other hand, Slodkowski's Lemma \BB.1  and Jensen's Lemma \GG.1 below are equivalent special cases of Theorem \AAA.7
(see Section \GG\  for a proof of this equivalence).

The four results above yield the following useful theorem concerning the
upper contact jets of a quasi-convex function.
The order two part of this theorem can 
be considered a ``partial upper semi-continuity of the
second derivative'' (PUSC of SD).

\Theorem{\AAA.8. (Upper Contact Jets)}
{\sl 
Suppose $u$ is quasi-convex with an upper contact jet $(p_0,A_0)$ at a point $x$.
 Then
\medskip
{\rm (D at UCP)}  \quad  $u$ is differentiable at $x$ and $D_{x}u = p_0$.
\medskip

Suppose $E$ is a set of full measure in a neighborhood of $x$. Then
there exists a sequence $\{x_j\} \ss E$ with $x_j\to x$  such that $u$ is twice differentiable
at each $x_j$ and 
\medskip
{\rm (PC of FD)}  \quad $D_{x_j} u \ \to\ D_{x} u=p_0$,
\medskip
 
{(\rm PUSC of SD)}  \quad $D_{x_j}^2 u \ \to\ A \leq A_0$.
}

\pf
By Alexandrov's Theorem, the set of points $x\in E$ where $u$ is twice differentiable,
is a set of full measure. In order to apply the Jensen-Slodkowski Lemma we replace 
$(p_0, A_0)$ by the strict upper contact jet $(p_0, A_0+\e I)$. Now choose a sequence
$\e_j\to0$, and pick a point $x_j \in B_{\e_j}(x_0)$ such that: (1) $x_j\in E$, (2) $u$ is twice differentiable at $x_j$, and (3) $x_j$ is a global upper contact point of type $A_0+\e_j I$
on $B_{\e_j}(x_0)$ for $u$.  By the basic differential calculus fact 
Lemma \AA.3, $D_{x_j}^2u \leq A_0+\e_j I$. 
Since $u$ is $\l$-quasi-convex, we have $D^2_{x_j} u +\l I \geq0$. Thus,
$$
-\l I \ \leq \ D^2_{x_j} u\ \leq\ A_0+\e_j I.
\eqno{(\AAA.4)}
$$
By compactness there is  a subsequence such that $D^2_{x_j} u \to A\leq A_0$.\qed
\medskip

Theorem  \AAA.8 can be stated succinctly in terms of the subset $J^+(u)\ss
J^2(X) \equiv X\times \bbr\times\rn\times \Symn$  of upper contact jets  for $u$, 
and another  subset depending on $E$.
Define  $J(u,E) \ss J^2(X)$  to be the subset of tuples $(x,u(x), D_x u, D^2_x u+P)$
 such that $x\in E$, $u$  is twice differentiable at $x$, and $P\geq0$.  Then Theorem \AAA.8  condenses to:
$$
{\rm If\ } u \ {\rm is \ quasi\!-\!convex\ and \ } E\ {\rm has\ full\  measure,\ then\ } \ J^+(u) \ \ss\ \overline{J(u,E)}.
\eqno{(\AAA.5)}
$$

   We will deduce the four  results  from the special case where $u$  is convex, and for Lemma 
   \AAA.7 we will reduce to the special case where $A=\l I$, i.e., Slodkowski's Lemma \BB.1.

%\vfill\eject
\vskip.3in

%%%%%%%%%%%%%%%%%%%%%%%%%%%%%%%%%%%%%%%%%%%%%%%%%%%%
%%%%%%%%%%%%%%%%%%%%%%%%%%%%%%%%%%%%%%%%%%%%%%%%%%%%
%%%%%%%%%%%%%%%%%%%%%%%%%%%%%%%%%%%%%%%%%%%%%%%%%%%%
%%%%%%%%%%%%%%%%%%%%%%%%%%%%%%%%%%%%%%%%%%%%%%%%%%%%
%%%%%%%%%%%%%%%%%%%%%%%%%%%%%%%%%%%%%%%%%%%%%%%%%%%%

\centerline{\headfont \BBB.\ Convex Functions -- The Subdifferential.}
 \medskip 

In this section  we shall assume that $u$ is convex and prove
some of the basic properties using the following standard concept involving lower contact.
If  a convex function $u$ is defined  on a convex open
set  $X\ss\rn$, the {\bf subdifferential $\partial u$ of $u$} is defined to be the set
$$
\partial u \ \equiv \ \{(x,p) \in X\times \rn : \ u(x)+\bra p {y-x} \leq u(y) \ \ \forall\, y\in X\}.
\eqno{(\BBB.1)}
$$
Geometrically this means that the graph of the affine function $u(x) +\bra p{y-x}$ lies below the graph of $u$ and the graphs touch above the point $x$, that is, the hyperplane is a supporting
hyperplane for the convex set $\{y\geq u(x)\}$ (or, in the language of Definition \AAA.5, that
$x$ is a lower contact point of type $A=0$ on $X$).

  The fibre of $\partial u$ over $x\in X$ is denoted
$\partial u(x)$.  Note that
$$
u\ \ {\rm is\ convex} \Iff \partial u(x) \ \neq\ \emptyset \ \ {\rm for\ each \ }\ x \in X.
\eqno{(\BBB.2)}
$$
For most of  the purposes of this section, this can be taken as the definition of convexity.

The useful inequalities
$$
\bra{p}{y-x}  \ \leq\ u(y)-u(x) \ \leq\ \bra q{y-x} \ \  \fa \ \  p\in \partial u(x),
q\in \partial u(y)
\eqno{(\BBB.3)}
$$
follow immediately from the definition (\BBB.1) of the set $\partial u$.  These inequalities,
stated geometrically, say that (with coordinates intrinsic to the affine line determined by
$x$ and $y$) if $x\leq y$, then $\bar p \leq s\leq \bar q$ where $s$ is the slope of the chord
above $[x,y]$ and $\bar p$ and $\bar q$ are the slopes of the supporting lines above 
$x$ and $y$.
 Note that (\BBB.3) implies
that $\partial u$ is monotone (non-decreasing).  That is,
$$
0\ \leq\ \bra{q-p}{y-x}\  \  \fa \ \  p\in \partial u(x),
q\in \partial u(y).
\eqno{(\BBB.4)}
$$

For each compact set $K\ss X$
$$
\partial u \cap (K\times \rn)\ \ {\rm is\ a\  nonempty\ compact\ set \  with\ convex\ fibres.}
\eqno{(\BBB.5)}
$$
\pf
Let $K_\d = \{x\in X: \dist(x,K) \leq \d\}$ and 
choose $\d$ small enough so that $K_\d  \ss X$.
Since $|u|$ is bounded on $K_\d$, the left hand inequality in (\BBB.3) gives the upper bound 
$$
\bra p{y-x} \ \leq\ 2\|u\|_{K_\d}\qquad \forall\, p\in\partial u(x), \ x\in K, \ {\rm and}\ y\in K_\d.
$$ 
Choosing  $y = x + \d {p\over |p|}$ gives 
$
\bra p {y-x} = \d|p| \leq 2\|u\|_{K_\d},
$
and so $|p|$ is bounded on $K$.  Since $\partial u \cap (K\times \rn)$ is closed and vertically
bounded, the compactness assertion follows.

Now fix $x$ and note that 
$$
\eqalign
{
\partial u (x) \ &=\ \{p : u(x) + \bra p {y-x} \leq u(y) \ \ \forall\, y\in X\}\cr
 &=\  \bigcap_{y\in X} \{p : L_y(p) \leq u(y) \}\cr
}
$$
is an intersection of affine half-spaces and  therefore convex.\qed

\medskip

Combining (\BBB.5) with the inequalities (\BBB.3) easily yields two important facts.
The first is that   $u$ is Lipschitz on $K$.

\Lemma{\BBB.1} {\sl 
$$
|u(y) -u(x)|\ \leq\ C|y-x|
\eqno{(\BBB.6)}
$$
where the Lipschitz constant $C$ is the supremum of $|p|$ taken over $p\in\partial u(x), x\in K$.}
\medskip

The second is the following.

\Lemma{\BBB.2} {\sl
$u$ is differentable at $x  \ \ \iff\ \  \partial u(x) =\{p\}$ is a singleton,
in which case}
$$
p\ =\ D_xu \ =\ \lim_{\matrix{y\to x\cr q\in \partial u(y)}} q.
\eqno{(\BBB.7)}
$$
\pf
If $\partial u(x) =\{p\}$ is a singleton, then by   the compactness in (\BBB.5)
$$
p  \ =\ \lim_{\matrix{y\to x\cr q\in \partial u(y)}} q.
\eqno{(\BBB.8)}
$$
With $p\in\partial u(x)$, $q\in\partial u(y)$ and $e=(y-x)/|y-x|$, the inequalities (\BBB.3)
can be rewritten as
$$
0\ \leq\ {u(y)-u(x) -\bra p{y-x}  \over |y-x|} \ \leq\ \bra {q-p} e.
\eqno{(\BBB.9)}
$$
Combining (\BBB.8) and (\BBB.9) shows that if $\partial u(x)$ is a singleton,
then $u$ is differentiable at $x$ with $D_xu=p$.

Suppose now that $u$ is differentiable at $x$ and 
$p\in \partial u(x)$.   Then $\bra p{y-x} \leq u(y)-u(x)$.
Hence, for each $e\in\rn, |e|=1$, we have $t\bra pe \leq u(x+te)-u(x)$ for $t$ small.
This proves that $\bra pe \leq \bra{D_xu} e$ for all $|e|=1$, and hence $p=D_xu$.
Since $\partial u(x)\neq \emptyset$, there always exists $p\in \partial u(x)$, 
and hence $D_xu\in\partial u(x)$.\qed

\vfill\eject

\Cor{\BBB.3}
\medskip
\centerline
{\sl
$u$ is differentiable everywhere $\ \ \iff \ \  \partial u$ is single valued  $\ \ \iff \ \  u$ is $C^1$.
}
\medskip

\Remark{ \BBB.4. (Critical Points)}
Note that by (\BBB.1)
\medskip
\centerline{
$x$ is a minimum point for $u
\qquad\iff\qquad
0\in \partial u(x)$
}
\medskip
\noindent
We say that $x$ is a {\sl critical point} for $u$ if $u$ is differentiable at $x$ and  $D_xu=0$.
By  Lemma \BBB.2
$$
 x \ {\rm is\  a \ critical\ point\ for\ } u \Iff
\partial u(x) \ =\ \{0\}.
\eqno{(\BBB.10)}
$$
In particular,
$$
 {\rm If }\ x \ {\rm is\  a \ critical\ point\ for\ } u, \ {\rm then}
 \ x \ {\rm is\  a \ minimum \ point\ for\ } u.
\eqno{(\BBB.11)}
$$
\medskip

Finally, in the proof of Alexandrov's Theorem it will be helpful to use additivity of the subdifferential
for the  sum of a convex function $u$ and a quadratic polynomial $\vf$ with $D^2\vf \geq0$.

\Lemma {\BBB.5} {\sl
Suppose $\vf(y) \equiv c+\bra qy +\half\bra{Py}y$ with $P\geq0$. Then}
$$
\partial(u+\vf)(x) \ =\ \partial u(x)+ \partial \vf(x)  \ =\ \partial u(x)  +q+Px.
$$
\pf  Of course, $\partial \vf (x) = \{q+Px\}$. Moreover, it
follows directly from the definition of the subdifferential
that $\partial u(x) +\partial \vf (x) \ss \partial(u+\vf)(x)$.
It remains the prove that $\partial(u+\vf)(x)  \subseteq\partial u(x)  +q+Px$.
Suppose $p+q+Px \in \partial (u+\vf)(x)$. We want to show that $p\in \partial u(x)$.
Our assumption is that  $u(y)+\vf(y) \geq u(x) +\vf(x) +\bra{p+q+Px}{y-x}$ which is equivalent
to 
\smallskip
\centerline
{
$u(y) \ \geq\ u(x) +\bra p{y-x} -\half \bra {P(y-x)}{y-x} \ \equiv\ \psi(y)$.
}
\smallskip
\noindent
This means that the epigraph   of $u$ is contained in the epigraph of $\psi$ 
in $\bbr^{n+1}$.  This remains true if one applies the dilation $\rho_t$ of $\bbr^{n+1}$ by
$t>0$ centered at the point $(x,u(x))$.  By convexity ${\rm epi}(u) \ss \rho_t({\rm epi}(u))$
for all $t\geq 1$.  That is,  ${\rm epi}(u) \ss \rho_t( {\rm epi}(u)) \ss \rho_t( {\rm epi}(\psi))$.
As $t\to\infty$, the dilations  $\rho_t({\rm epi}(\psi))$ decrease down to the half-space 
$H\equiv  \{y: u(x) +\bra p{y-x}\geq0\}$, proving that ${\rm epi}(u) \subseteq H$ or that $p\in \partial u(x)$.\qed

%\vfill\eject
\vskip.3in

%%%%%%%%%%%%%%%%%%%%%%%%%%%%%%%%%%%%%%%%%%%%%%%%%%%%
%%%%%%%%%%%%%%%%%%%%%%%%%%%%%%%%%%%%%%%%%%%%%%%%%%%%
%%%%%%%%%%%%%%%%%%%%%%%%%%%%%%%%%%%%%%%%%%%%%%%%%%%%
%%%%%%%%%%%%%%%%%%%%%%%%%%%%%%%%%%%%%%%%%%%%%%%%%%%%
%%%%%%%%%%%%%%%%%%%%%%%%%%%%%%%%%%%%%%%%%%%%%%%%%%%%

\centerline{\headfont \AA.\ Proof of Lemmas \AAA.2 and \AAA.3.}
 \medskip 

These two results follow easily from the convex case, enabling us to assume that 
 $u$ is convex in the proofs.
 
\medskip
\noindent
{\bf Proof of Lemma \AAA.2  (D at UCP).}
By (\BBB.2) there exists $\bar p\in\rn$ with $u(x) + \bra { \bar p} {y-x} \leq u(y), \ \forall y\in X$.
Subtracting the affine function of $y$ on the left from $u$, we can assume
$0\leq u(y)$ and $u(x)=0$.  
Now if $(p,A)$ is an upper contact jet for $u$ at $x$, then
$$
 0\ \leq\ u(y) \ \leq\   \bra p {y-x}  + \half \bra{A(y-x)}{y-x}
\eqno{(\AA.1)}
$$
which implies that $u$ is differentiable at $x$ with $D_xu =p$. This proves Lemma \AAA.2.\qed

\Remark{\AA.1}
This argument proves more.  First note that by (\BBB.1) if $0\leq u(y)$ and $u(x)=0$, then
$0\in \partial u(x)$. Now because of (\BBB.10) the inequalities (\AA.1) imply that $p=0$ and $A\geq0$.  Therefore,
$$
{\rm If } \ (p,A) \ {\rm is \  an \ upper\ contact\ jet\ for\ a \ convex\ function\ } u,\ {\rm\  then \ } A\geq0.
\eqno{(\AA.2)}
$$
For the converse, namely:
\medskip
\centerline
{
If each upper contact jet $(p,A)$ of an u.s.c. function $u$ satisfies $A\geq0$, then $u$ is convex,
}
\medskip
\noindent
the reader is referred to [HL$_1$].

\medskip
\noindent
{\bf Proof of Lemma \AAA.3 (PC of FD).}
By Lemma \BBB.2
$$
 {\rm If }\ u \ {\rm is\ differentiable \ at} \ x, \ {\rm then} \ \ 
 \lim_{\eqalign{\ \ \ \ &\ \ y\to x \ \cr&q\in \partial u (y)}} q\ =\ D_xu.
 \eqno{(\AA.3)}
$$
This is, in fact, a stronger version of Lemma \AAA.3.\qed

\vfill\eject
%\vskip .3in

%%%%%%%%%%%%%%%%%%%%%%%%%%%%%%%%%%%%%%%%%%%%%%%%%%%%
%%%%%%%%%%%%%%%%%%%%%%%%%%%%%%%%%%%%%%%%%%%%%%%%%%%%
%%%%%%%%%%%%%%%%%%%%%%%%%%%%%%%%%%%%%%%%%%%%%%%%%%%%
%%%%%%%%%%%%%%%%%%%%%%%%%%%%%%%%%%%%%%%%%%%%%%%%%%%%
%%%%%%%%%%%%%%%%%%%%%%%%%%%%%%%%%%%%%%%%%%%%%%%%%%%%

\centerline{\headfont \BB.\ The Reduction of the Jensen-Slodkowski Lemma to the Slodkowski Lemma.}
 \bigskip 

The Jensen-Slodkowski Lemma \AAA.7 contains the following as a special case.

\Lemma{\BB.1. (Slodkowski)} {\sl
Suppose that $u$ is a convex function with a strict upper contact jet $(0,\l I)$
at a point $x$.  Then there exists $\bar \rho >0$ such that the measure}
$$
\left| \cc(u, B_{\rho}, \l I)    \right| \ >\ 0 \qquad\forall\, 0<\rho\leq\bar \rho.
$$

The following trivial lemma is all that is needed for the reduction.
First note that 
a degree-2 polynomial $\vf(y)$ satisfies
$$
\vf(y)\ =\ \vf(x) + \bra{D_x\vf}{y-x} +\half \bra{(D^2_x \vf) (y-x)}{y-x} 
\qquad\forall\, x,y\in\rn.
$$
and $D^2_x\vf$ is independent of $x$.

\Lemma{\BB.2} {\sl
Suppose $\vf$ is a degree-2 polynomial.  Set $B\equiv D^2_x\vf$.
\medskip
(1) \ \ If $(p,A)$ is an upper contact jet for $u$ at $x$ on $X$, then
$(p+D_x\vf, A+B)$ is an upper contact jet for $w\equiv u+\vf$ at $x$ on $X$.

\medskip
(2) \ \ $(p,A)$ is strict for $u$ \ \ $\Rightarrow$\ \ $(p+D_x\vf, A+B)$ is strict for $w\equiv u+\vf$.

\medskip
(3) \ \ $\cc(u+\vf, X, A+B) \ =\ \cc(u,X,A)$
}
\pf
For any point $x\in X$
$$
\eqalign
{
u(y)\ &\leq \ u(x) + \bra{p}{y-x} + \half \bra{A(y-x)}{y-x} \quad \forall\, y\in X  \cr
\iff\qquad
w(y)\ &\leq \ w(x) + \bra{p+D_x\vf}{y-x} + \half \bra{(A+B)(y-x)}{y-x} \quad \forall\, y\in X 
}
$$
\hfill\qed

\medskip
We now claim the following.
$$
{\sl The \ special\ case \ Lemma\ \BB.1\ implies\ the\ full\  Jensen-Slodkowski \ Lemma.}
\eqno{(\BB.3)}
$$
\noindent
{\bf Proof.}
Suppose $(p,A)$ is a strict upper contact jet for $u$ at $x$.
Take $\vf(y) \equiv  -\bra {p}{y-x} - \half \bra{A(y-x)}{y-x} + {\l\over 2} |y-x|^2$
and apply Lemma \BB.2.  Then $(0,\l I)$ is a strict upper contact jet for $w\equiv u+\vf$
at $x$ on $X \equiv B_{\bar\rho}(x)$ for some $\bar \rho>0$.
Moreover,  $\cc(w,X,\l I) = \cc(u,X,A)$.
Finally take $\l$ sufficiently large so that $w\equiv u+\vf$ is convex.
(If $u$ is $\a$-quasi-convex and $A\leq \b I$, take $\l\geq \a+\b$.)
Now  the Slodkowski Lemma \BB.1 can be applied to $w$.\qed
\medskip

In the next section we establish the elementary convex geometric fact needed to prove
Slodkowski's Lemma.

\vfill\eject
%\vskip.3in
%\vfill\eject

%%%%%%%%%%%%%%%%%%%%%%%%%%%%%%%%%%%%%%%%%%%%%%%%%%%%
%%%%%%%%%%%%%%%%%%%%%%%%%%%%%%%%%%%%%%%%%%%%%%%%%%%%
%%%%%%%%%%%%%%%%%%%%%%%%%%%%%%%%%%%%%%%%%%%%%%%%%%%%
%%%%%%%%%%%%%%%%%%%%%%%%%%%%%%%%%%%%%%%%%%%%%%%%%%%%
%%%%%%%%%%%%%%%%%%%%%%%%%%%%%%%%%%%%%%%%%%%%%%%%%%%%

\def\epi{{\rm epi}}
\def\CH{{\rm ch}}
\def\SLAB{{\rm SLAB}}
\def\gr{{\rm graph}}

\centerline{\headfont \CC.\ The Convex Hull of Two Open Paraboloids of the Same Radius.}
 \medskip

Now we begin the proof of Slodkowski's Lemma \BB.1, which is completed in Section \FF.
Our proof is an adaptation of his proof, in which we use  paraboloids of radius $r$
in place of balls of radius $r$.  
It is important that these paraboloids be {\sl open}. Each such paraboloid is determined
by its  vertex $(v, \vf(v))\in \rn\times \bbr$, and by definition, is the open epigraph
$\epi(\vf)$ of the quadratic function
$$
\vf(y)\ \equiv\ \vf(v) + {1\over 2r}|y-v|^2.
$$
Given two such open paraboloids $\epi(\vf_1)$ and $\epi(\vf_2)$
with vertices $(v_1,\vf(v_1))$ and  $(v_2,\vf(v_2))$ respectively, we compute the convex 
hull $\CH(\epi(\vf_1) \cup \epi(\vf_2))$
of the union of these two open sets.  We shall emphasize what is needed in the application.

\Lemma{\CC.1} {\sl
There is an open vertical slab $\SLAB\ss\bbr^{n+1}$ written as the intersection 
$\SLAB = \ch_1 \cap \ch_2$ of two parallel vertical open half-spaces with the following
property.  Let $\CH \equiv \CH(\epi(\vf_1) \cup \epi(\vf_2))$.
Then
$$
\gr(\vf_1) \cap \CH \ \ss\ \ch_1
\and
\gr(\vf_2) \cap \CH \ \ss\ \ch_2
\eqno{(\CC.1)}
$$
Moreover, the width of $\SLAB$ is $|v_1-v_2|$.
}
\pf
Set $e\equiv {v_2-v_1\over |v_2-v_1|}$ and let $m\equiv {\vf(v_2)- \vf_1(v_1)\over |v_2-v_1|}$
denote the slope of the line segment from the first vertex to the second.
Define  $\ch_1$ to be  the open half-space whose boundary hyperplane
$\partial \ch_1$ has   interior normal $(e,0)$  and passes through $(v_1+r me, 0)$.
Similarly, define $\ch_2$ to have  interior  normal $(-e,0)$ and boundary $\partial \ch_2$
passing through $(v_2+r me, 0)$.  Then $\SLAB \equiv \ch_1\cap \ch_2$ clearly  has 
width  $|v_2-v_1|$.   It remains to prove (\CC.1).

This can be seen  by determining the pairs of points 
$$
z_1 \ \equiv \ (y_1, \vf_1(y_1))\in \gr(\vf_1)
\and
z_2 \ \equiv \ (y_2, \vf_2(y_2))\in \gr(\vf_2)
$$
which have a common tangent plane $H$.  Equating normals $(D_{y_1} \vf_1, -1)$ and
$(D_{y_2} \vf_2, -1)$ yields $y_1-v_1=y_2-v_2$.  Thus $y_1=v_1+w$ and $y_2=v_2+w$
for some $w\in \rn$.

Now $N\equiv ({w\over r}, -1)$ is normal to $H$.  Hence,  $z_1, z_2 \in H$ implies that
$
{1\over r} \bra {y_1}{w} - \vf_1(y_1) =  {1\over r} \bra {y_2}{w} - \vf_2(y_2) 
$.
Therefore $\bra {v_2-v_1}w = \bra {y_2-y_1}w  = r(\vf_2(y_2) -\vf_1(y_1))$.
However,  
$$
\vf_2(y_2) -\vf_1(y_1) = \vf_2(v_2)  + {1\over 2r}|y_2-v_2|^2 -\vf_1(v_1) -  {1\over 2r}|y_1-v_1|^2
= \vf_2(v_2)-\vf_1(v_1)
$$
proving that $\bra ew =rm$. Let $\bbr^{n-1}$ denote $e^\perp$ in $\rn$. 
This proves that there exists   $\bar w\in\bbr^{n-1}$ with
$$
z_1\ =\ (v_1+rme+\bar w, \vf_1(v_1+rme+\bar w))
\and
z_2\ =\ (v_2+rme+\bar w, \vf_1(v_2+rme+\bar w)).
$$
The mapping $\bar w\to z_1$ with $\bar w\in \bbr^{n-1}$ parameterizes $\partial \ch_1\cap \gr(\vf_1)$,
and similarly  $\bar w\to z_2$  parameterizes $\partial \ch_2\cap \gr(\vf_2)$.  \qed

\Remark{\CC.2}  Consider the closure $C$ of $\CH(\epi(\vf_1) \cup \epi(\vf_2))$.
The points in $\epi(\vf_1) \cup \epi(\vf_2) \sim \SLAB$ are extreme points of $C$.
For each $\bar w\in\bbr^{n-1}$ as above, the associated hyperplane $H$
supports $C$ and intersects $C$ along the line segment from $z_1$ to $z_2$.

\vfill\eject
%\vskip .3in

%%%%%%%%%%%%%%%%%%%%%%%%%%%%%%%%%%%%%%%%%%%%%%%%%%%%
%%%%%%%%%%%%%%%%%%%%%%%%%%%%%%%%%%%%%%%%%%%%%%%%%%%%
%%%%%%%%%%%%%%%%%%%%%%%%%%%%%%%%%%%%%%%%%%%%%%%%%%%%
%%%%%%%%%%%%%%%%%%%%%%%%%%%%%%%%%%%%%%%%%%%%%%%%%%%%
%%%%%%%%%%%%%%%%%%%%%%%%%%%%%%%%%%%%%%%%%%%%%%%%%%%%

\centerline{\headfont \DD.\  Upper Semi-Continuous Functions -- Radius-r Upper Contact Points.}
 \medskip 

Analogous to the fact that  the subdifferential is basic for understanding convex functions,
is the fact that upper contact quadratics of radius $r$ are basic for understanding general 
upper semi-continuous functions.  In this section $u$ is any 
upper semi-continuous function on $X$.
Let $\vf$ be a quadratic function, and note that its second derivative $D^2_x\vf$
is independent of the point $x$.  We say that $\vf$ has {\bf radius $r$} if $D^2_x\vf = {1\over r}I$
where $0<r<\infty$.  In this case $\vf$ has a unique minimum point $v$ which we call the
{\bf vertex point}.  We say the graph of $\vf$ has its {\bf vertex} at $(v, \vf(v))$.
The radius $r$ and vertex $v$ determine $\vf$ up to its {\bf height} $\vf(v) =c$, that is
$$
\vf(y)\ =\ c + {1\over 2r} |y-v|^2
\eqno{(\DD.1)}
$$

Now assume that $X$ is a compact set in $\rn$ and $u$ is an 
upper semi-continuous function on $X$.  For large $c$ the graph of $\vf$
lies above the graph of $u$.  The smallest such $c$ is
$$
\eqalign
{
\wh c \ &\equiv\ \inf \left \{ c : u(y) \leq c + \smfrac 1 {2r} |y-v|^2 \ \forall\, y\in X \right \}  \cr
&=\ \sup_{y\in X} \left( u(y) - \smfrac 1 {2r} |y-v|^2  \right).
}
\eqno{(\DD.2)}
$$
Since $u$ is upper semi-continuous and $X$ is compact, this supremum
$\wh c$ is attained at some point $x\in X$.  Thus with the height $c=\wh c$
defined by (\DD.2), the polynomial $\vf(y)$ given by (\DD.1) satisfies
$$
\eqalign
{
 &{\rm(a)} \ \ u(y)\ \leq\ \vf(y) \quad {\rm for\ all\ } y\in X, \ {\rm and}    \cr
&{\rm(b)} \ \ u(x)\ = \vf(x) \quad {\rm for\ some\ } x\in X. 
}
\eqno{(\DD.3)}
$$
The contact set $\{u=\vf\}$ is compact since $X$ is compact.  It will be denoted by
$$
\cc(u, X, \smfrac 1r I, v).
\eqno{(\DD.4)}
$$

If $v$ belongs to the interior of $X$ and  the radius $r$ is small, these contact points $x$, where graph$(\vf)$ touches
graph$(u)$, should occur in the interior of $X$.  The basic estimate proving this is given in the following lemma.  Let   ${\rm Osc}_X(u) \equiv \sup_X u-\inf_X u$ denote the oscillation of $u$ on $X$.
Note that ${\rm Osc}_X(u) <\infty$ if and only if $u$ is bounded below, since $u$ is 
upper semi-continuous.  Finite oscillation must be assumed in order for the estimate to have
content.

\Lemma{\DD.1} {\sl
If $x\in \cc(u, X, \smfrac 1r I, v)$, then
}
$$
|x-v|  \ \leq\ \sqrt{2r {\rm Osc}_X(u)}.
$$
%\pf
%We have $\vf(y) \equiv \vf(v) + {1\over 2r} |y-v|^2 \geq u(y)$ for all $y$, and $\vf(x)=u(x)$.
%Thus, $ {1\over 2r} |x-v|^2  = \vf(x)-\vf(v) = u(x) -\vf(v) \leq u(x) -u(v) \leq {\rm Osc}_X(u)$.\qed
\pf
 Since the graph of $\vf$ lies above the graph of $u$ and $\vf(x) = u(x)$,
the change in $\vf$ from $v$ to $x$, which equals $\vf(x) -\vf(v) = {1\over 2r}|x-v|^2$,
is $\leq$ the change   $u(x)-u(v)$,  which is $\leq$ Osc$_X(u)$.\qed
\medskip

This result can be put in a more useful form, guaranteeing lots of upper contact points
for any upper semi-continuous function which is bounded below.

\Lemma{\DD.1$'$} {\sl
Set $\d \equiv \sqrt{2r{\rm Osc}_X(u)}$ and $X_\d \equiv \{y\in X : \dist(y, \partial X)>\d\}$.
For  any point $v\in X_\d$ the contact set $\cc(u, X, \smfrac 1r I, v)$ is a non-empty compact subset 
of the open set $X_\d$.  In fact, it is contained in the closed ball $B_\d(v)$ about $v$ of radius $\d$.
}
\bigskip

\centerline{\bf The Upper Vertex Map}
\medskip

Now suppose $X\ss\rn$ is open and 
$$
u(y) \ \leq\ u(x) + \bra p{y-x} +{1\over 2r}|y-x|^2  \quad \forall\, y\in X.
\eqno{(\DD.5)}
$$
By Definition \AAA.5 $x \in \cc(u,X, \smfrac 1r I)$, that is, $x$ is a global upper contact point 
of type ${1\over r}I$ (radius $r$)
on $X$ for $u$.  We also say that $(p, {1\over r} I)$ is an upper contact jet for $u$ at $x$
with the upper contact inequality holding on all of $X$.  

Set
$$
\vf(y) \ \equiv\ u(x) + \bra p{y-x} +{1\over 2r}|y-x|^2.
\eqno{(\DD.6)}
$$
Then $\vf$ has radius $r$, and the vertex point $v$ for $\vf$ is given by
$$
v\ =\ x-rp,
\eqno{(\DD.7)}
$$
since $0=D_v\vf = p+{1\over r}(v-x)$. Furthermore, if $x\in {\rm Diff}^1(u)$, the set of points
where $u$ is differentiable, 
then the $p$ satisfying (\DD.5) is unique and equal to $D_xu$.

\Def{\DD.3. (The Upper Vertex Map)}   The map
$$
V:\cc (X,u,\smfrac 1{r} I)\cap {\rm Diff}^1(u) \ \arr\ \rn
$$
defined by $V(x) \equiv x-rD_x u$, i.e., $V\equiv I-rDu$, will be called the 
{\bf vertex map for} $u$.  

This map has the property that
$$
v\ =\ V(x) \qquad\Rightarrow\qquad |v-x|\ \leq\ \sqrt{2r{\rm Osc}_X(u)}
\eqno{(\DD.8)}
$$
by the basic estimate Lemma \DD.1.

\bigskip
\centerline
{
\bf Convex Functions
}
\medskip

The constructions of the previous section apply to a convex function $u$.
By (D at UCP) we have $\cc(u,X, {1\over r}I) \ss {\rm Diff}^1(u)$.

Therefore, the vertex map $V \equiv I-r Du$ is a well defined map
$$
V : \cc(u, X, \smfrac 1r I) \ \arr\ \rn.
\eqno{(\DD.9)}
$$

\vfill\eject
%\vskip .3in

%%%%%%%%%%%%%%%%%%%%%%%%%%%%%%%%%%%%%%%%%%%%%%%%%%%%
%%%%%%%%%%%%%%%%%%%%%%%%%%%%%%%%%%%%%%%%%%%%%%%%%%%%
%%%%%%%%%%%%%%%%%%%%%%%%%%%%%%%%%%%%%%%%%%%%%%%%%%%%
%%%%%%%%%%%%%%%%%%%%%%%%%%%%%%%%%%%%%%%%%%%%%%%%%%%%
%%%%%%%%%%%%%%%%%%%%%%%%%%%%%%%%%%%%%%%%%%%%%%%%%%%%

\centerline{\headfont \EE.\  The Vertex Map is a Contraction.}
 \medskip

\Prop{\EE.1} {\sl
Given a convex function $u$  defined on an open convex set $X\ss\rn$,
the vertex map $V : \cc(u,X, {1\over r} I))\to\rn$ is a contraction.
}

\Remark{\EE.2}  If $u$ is smooth, then the Jacobian $J$ of $V$ is $I-rD^2 u$
and $0\leq J\leq I$ on $\cc(u,X, {1\over r} I))$ is obvious.

\pf Given $x_1, x_2\in \cc(u,X, {1\over r} I))$ we must show that
$$
|V(x_2)-V(x_1)|\ \leq\ |x_2-x_1|.
\eqno{(\EE.1)}
$$
Let $\vf_1(y)$ denote the quadratic of radius $r$ whose graph lies above $\gr(u)$
and touches at $(x_1, u(x_1))$, i.e., $\vf_1(x_1)=u(x_1)$, and define $\vf_2$ similarly.
Let $v_1\equiv  V(x_1)$ and $v_2 \equiv V(x_2)$ denote the vertex points.  

Since $u$ is convex and $\vf_k\geq u$, $k=1,2$, we have 
$$
 \CH\left( \epi(\vf_1) \cup \epi(\vf_2)\right) \ \ss\ \epi(u).
$$
Now $(x_1, \vf_1(x_1)) = (x_1, u(x_1)) \notin \epi(u)$ (recall that $\epi(u)$
is the {\sl open} epigraph).   Hence, 
$(x_1, \vf_1(x_1)) \notin  \CH\left( \epi(\vf_1) \cup \epi(\vf_2)\right)$,
and therefore $(x_1, \vf_1(x_1))\notin \ch_1$ by Lemma \CC.1.
Similarly, $(x_2, \vf_1(x_2)) \notin \ch_2$. We conclude that these points lie on opposite sides
of $\SLAB$ and so $|x_2-x_1|\geq {\rm width}(\SLAB) = |v_2-v_1|$.\qed

%\vfill\eject
\vskip .3in

%%%%%%%%%%%%%%%%%%%%%%%%%%%%%%%%%%%%%%%%%%%%%%%%%%%%
%%%%%%%%%%%%%%%%%%%%%%%%%%%%%%%%%%%%%%%%%%%%%%%%%%%%
%%%%%%%%%%%%%%%%%%%%%%%%%%%%%%%%%%%%%%%%%%%%%%%%%%%%
%%%%%%%%%%%%%%%%%%%%%%%%%%%%%%%%%%%%%%%%%%%%%%%%%%%%
%%%%%%%%%%%%%%%%%%%%%%%%%%%%%%%%%%%%%%%%%%%%%%%%%%%%

\centerline{\headfont \FF.\  Completion of the Proof of Slodkowski's Lemma.}
 \medskip 
 
 \def\Osc{\rm Osc}
 
 The fact that the vertex map is a contraction combined with a standard perturbation
 argument is all that is needed to prove Slodkowski's Lemma  \BB.1.
 
 We may assume that $x=0$ and $u(x)=0$.  Then the assumption on the convex
 function $u$ which occurs in Slodkowski's Lemma  is:
 $$
 (0,\l I)\ \ {\rm is\ a \ strict\ upper\ contact \ jet\ for\ } u \ {\rm at  \ } x.
 \eqno{(\FF.1)}
 $$
With  $\l=1/R$ we claim that this is equivalent to the existence of $\bar \rho>0$ such that:
$$
 0\ \leq\  u(y) \ <\ {1\over 2R} |y|^2\qquad {\rm for\ } 0<|y|\leq \bar \rho.
 \eqno{(\FF.1)'}
 $$
Assuming (\FF.1), the function $u$ is differentiable at $x_0=0$ and $D_{x_0}u=0$
by (D at UCP).  The convexity of $u$ then implies that $0\leq u(y)$ (cf.(\BBB.5)), showing that 
(\FF.1) $\Rightarrow$ (\FF.1)$'$. That (\FF.1)$'$ $\Rightarrow$ (\FF.1) is clear.

Now we apply Lemma \DD.1$'$ to $u$ with $X\equiv \overline{B}_\rho$ 
where $\rho\leq \overline\rho$ and $r\leq R$.
Since $u(0) = 0$ and $u(y) \geq 0$ on  $\overline{B}_\rho$, we have 
$\Osc_X = \sup_{\overline{B}_\rho} u$, which we denote by $M(\rho)$
so that $\d(\rho) \equiv \sqrt{2rM(\rho)}$.  Now $X(\d) = B_{\rho-\d}$.
According to Lemma \DD.1$'$, if $v\in  B_{\rho-\d}$, then $\cc (u, B_\rho, {1\over r} I, v)$ 
is a non-empty subset of $B_{\rho-\d}$.  Choose such an upper contact point $x$.  
Then, as noted above, since $u$ is convex, the vertex map $V$ is defined at $x$.  
Thus for all $v\in B_{\rho-\d}$, we have that 
$v=V(x)$  for some $x\in \cc (u, \overline B_\rho, {1\over r} I)$,
i.e.,   $\overline B_{\rho-\d} \ss V(\cc (u, \overline B_\rho, {1\over r} I))$.
Since $V$ is a contraction, this proves that: 
$$
|\overline B_{\rho-\d}|\ \leq \ |V(\cc (u, \overline B_\rho,  \smfrac 1 r  I))| \ \leq \ | \cc (u, \overline B_\rho,  \smfrac 1 r  I) |.
\eqno{(\FF.2)}
$$

Since $M(\rho) = \sup_{\partial B_\rho} u$ by the maximum principle,
(\FF.1)$'$  implies that $M(\rho) < \rho^2/(2R)$, and hence, $\d(\rho) < \rho \sqrt{r\over R}$.
This proves that 
$$
\rho - \d(\rho) \ >\ \rho \left(  1-\sqrt{ \smfrac rR}\right).
\eqno{(\FF.3)}
$$
Taking $r=R$, this completes the proof of Lemma \BB.1, since $\rho-\d(\rho)>0$.\qed
\medskip

Our proof of this lemma, which uses paraboloids instead of spheres, also provides an analogue of
Slodkowski's Lemma \AA.4 which estimates the density of the contact set,
using (\FF.3).

\Prop{\FF.1}  {\sl
If $u$ is a convex function on $B_{\overline \rho}$ satisfying 
$0\leq u(y) < {|y|^2\over 2R}$ for $y\neq 0$.  Then for $0< r< R$, }
$$
|B_{\rho} | \left(1-\sqrt{\smfrac r R} \right)^n \ \leq \ \left | \cc \left(u, \overline B_{\rho}, \smfrac 1 r I \right)\right |
\qquad \forall\, 0<\rho \leq \overline \rho.
 $$

%\vfill\eject
\vskip .3in

%%%%%%%%%%%%%%%%%%%%%%%%%%%%%%%%%%%%%%%%%%%%%%%%%%%%
%%%%%%%%%%%%%%%%%%%%%%%%%%%%%%%%%%%%%%%%%%%%%%%%%%%%
%%%%%%%%%%%%%%%%%%%%%%%%%%%%%%%%%%%%%%%%%%%%%%%%%%%%
%%%%%%%%%%%%%%%%%%%%%%%%%%%%%%%%%%%%%%%%%%%%%%%%%%%%
%%%%%%%%%%%%%%%%%%%%%%%%%%%%%%%%%%%%%%%%%%%%%%%%%%%%

\centerline{\headfont \GG.\  The Equivalence of Slodkowski's Lemma and Jensen's Lemma.}
 \medskip

The results of this section are not needed elsewhere in these notes, but they might have some 
historical interest.

Another special case of the general Slodkowski-Jensen Lemma \AAA.7 is Jensen's Lemma.

\Lemma{\GG.1.  (Jensen)} {\sl
Suppose that $w$ is a quasi-convex function with the strict  upper contact jet $(0,0)$ at $x$
(equivalently, $w$ has a strict local maximum at $x$).
Then there exists $\bar \rho>0$ such that 
$$
\left| \cc(w, B_{\rho}, 0)    \right| \ >\ 0 \qquad\forall\, 0<\rho\leq\bar \rho.
$$
}
$${\rm
 In \ a\  very\  strong \ sense\  the\  Slodkowski\  Lemma \  is \ equivalent\  to\  
Jensen's\  Lemma.}
\eqno{(\GG.1)}
$$
The precise statement (\GG.1) is embedded in the following proof.
\pf
Set
$$
u(y)\ \equiv\ w(y) + {\l\over 2} |y-x|^2.
\eqno{(\GG.2)}
$$
Then by definition
\medskip
\centerline
{
$u$ is convex \qquad$\iff$\qquad $w$ is $\l$-quasi-convex.
}
\medskip

By Lemma \BB.2 parts (1) and (2)
$$
\eqalign
{
&(0,\l I) \ \ {\rm is\ a\ strict\ upper\ contact\ jet\ for\ } u \ {\rm at \ } x \qquad\iff\qquad \cr
&(0,0) \ \  {\rm is\ a\ strict\ upper\ contact\ jet\ for\ } w\  {\rm at \ } x,
}
\eqno{(\GG.3)}
$$
 while part (3) state that
$$
\cc(u, B_\rho, \l I) \ =\ \cc(w, B_\rho, 0).
\eqno{(\GG.4)}
$$
Thus by (\GG.3) the hypotheses of Slodkowski and Jensen are equivalent, while by
(\GG.4) the conclusions of Slodkowski and Jensen are identical (not just equivalent).

%\vfill\eject
\vskip .3in

%%%%%%%%%%%%%%%%%%%%%%%%%%%%%%%%%%%%%%%%%%%%%%%%%%%%
%%%%%%%%%%%%%%%%%%%%%%%%%%%%%%%%%%%%%%%%%%%%%%%%%%%%
%%%%%%%%%%%%%%%%%%%%%%%%%%%%%%%%%%%%%%%%%%%%%%%%%%%%
%%%%%%%%%%%%%%%%%%%%%%%%%%%%%%%%%%%%%%%%%%%%%%%%%%%%
%%%%%%%%%%%%%%%%%%%%%%%%%%%%%%%%%%%%%%%%%%%%%%%%%%%%

\centerline{\headfont \HH.\  The Proof of Alexandrov's Theorem.}
 \medskip 
 
 We will evoke two local results about Lipschitz maps $G:\rn\to\rn$
 which can be found many places (e.g. [F]).
 Otherwise the proof is elementary and complete.
 It combines elements of the proofs in [CIL] and [AA] with the Legendre transform.
 It is also worth noting that while upper contact quadratics of radius $r$ were key to the 
 proof of Slodkowski's Lemma, lower contact points of radius $r$ are the key to proving
 Alexandrov's Theorem.

 \medskip
 \noindent
 {\bf  Rademacher's Theorem.}
{\sl
The derivative of $G$ exists almost everywhere.
}
\medskip

 It is convenient to lable the variables as $x=G(y)$.
We will say  that $y$ {\sl is a  critical point for} $G$ if $G$ is differentiable at $y$ and $D_yG$
is singular.  The image  of the set of all critical points under the mapping $G$
is the {\sl set of critical values of $G$}.

 \medskip
 \noindent
 {\bf  The Lipschitz Version of Sard's Theorem.}
{\sl
The set of critical values of $G$ has measure zero.
}
\medskip

Since convex functions are locally Lipschitz (see Lemma \BBB.1), the scalar version of Rademacher's Theorem implies that
$$
{\rm A\ convex\ function\ is\ differentiable \  almost \  everywhere.}
\eqno{(\HH.1)}
$$
This fact will also be used in the proof.

Now we begin the proof of Alexandrov's Theorem.  Given a convex function $u$
on a convex open set $X$ with upper bound $N$ and $r>0$,
we will show that 
$$
f(x) \ \equiv \  r u(x)  + \half |x|^2
$$
 is twice differentiable a.e.
on $X_\d \equiv \{x\in X : \dist(x,\partial X)>\d\}$ where $\d\equiv 4\sqrt{r|u|_{\infty}}$
and $|u|_{\infty} \equiv \sup_X |u|$.  
First note that by Lemma \BBB.5
$$
\partial f(x) \ =\ x+r\partial u(x), \ \ \  {\rm   i.e.,}\ \ \  \partial f\ =\ I+r\partial u.
\eqno{(\HH.2)}
$$
That is, with $x\in X$,
$$
{\rm If \ } y\ {\rm and}\ p \ {\rm are\  related\ by \ } y=x+rp,\quad {\rm then\ } \ \ 
(x,y)\in \partial f
\ \ \iff\ \
(x,p)\in \partial u.
\eqno{(\HH.3)}
$$
\Lemma{\HH.1} {\sl
The multi-valued map $\partial f$ is expansive. That is,
if $(x_1,y_1)\in \partial f$ and $(x_2,y_2)\in \partial f$, then
}
$$
|x_1-x_2|\ \leq |y_1-y_2|.
$$
\pf
Note that $ |y_1-y_2| |x_1-x_2| \geq \bra {y_1-y_2}{x_1-x_2} =
|x_1-x_2|^2 + r \bra{p_1-p_2}{x_1-x_2}$
with $p_1$ and $p_2$ defined by (\HH.3) so that 
$p_1\in \partial u(x_1)$ and $p_2\in \partial u(x_2)$.  
By the monotonicity (\BBB.4) of $\partial u$ this expression is $\geq |x_1-x_2|^2$.\qed

\medskip

Because of this inequality, if $(x_1,y), (x_2,y) \in \partial f$, then $x_1=x_2$.
Thus the inverse of $F \equiv \partial f$ is single-valued.
We denote this single-valued mapping, which is defined on the set
\medskip
\centerline
{
$Y \ \equiv\ {\rm Im} F \ \equiv\ $ the projection
of $\partial f$ onto the second factor of $\rn\times\rn$,}
\medskip
\noindent
 by $x=G(y)$.
Now Lemma \HH.1 states that
$$
 \ |G(y_1)-G(y_2)|\ \leq\ |y_1-y_2|  \qquad{\rm for}\ \ y_1,y_2\in Y.
\eqno{(\HH.4)}
$$
That is, $G$ {\sl is  1-Lipschitz}  or {\sl contractive}. % (with $\leq$ not $<$).

\Lemma{\HH.2} 
$$
X_\d\ \ss\ {\rm Dom(G)} \ \equiv\  Y\ \ \ {\rm where}\ \ \d=2\sqrt{r|u|_{\infty}}.
\eqno{(\HH.5)}
$$
\pf Given $y\in X_\d$ pick a minimum point $x$ for the function $ru(z)  + \half |z-y|^2$
on the ball $B_\d(y) \equiv \{x : |x-y|\leq\d\} \ss X$.
This is also a minimum for the same function on all of $X$.
To see this consider any $\bar x \in X$ with  $|y-\bar x|>\d$. 
Then $u(\bar x)-u(y) +{1\over 2r}|\bar x-y|^2 > -2|u|_\infty + {\d^2\over 2r} =0$
by the definition of $\d$. This implies that  the value of $ru(z)+\half |z-y|^2$ at $\bar x$
is greater  than the value at $y$, and therefore greater than the value at $x$.  

The fact that  $x$ is a minimum point of $ru(z) +\half|z-y|^2$ on $X$
implies that $0\in  \partial(ru(z)+\half|z-y|^2)(x)$.  By Lemma \BBB.5 there 
exists $p\in \partial u(x)$ with $0= rp+(x-y)$. By  (\HH.3) above, $y\in \partial  f(x)$.
\qed

\bigskip
\centerline{\bf The Legendre Transform}\medskip

The map $G$ is the inverse of the multi-valued map $F=\partial f$. This map $G$ also has
a scalar potential $g$, with $G=\partial G$, which is classically called the Legendre transform
of $f$. Let $Y={\rm Im} F$ as above an note that with $(x,y)\in X\times Y$
$$
\eqalign
{
y\in\partial f(x) \Iff &f(x) + \bra y {z-x} \leq f(z) \ \ \ \forall \, z\in X  \cr
\Iff &f(x) -  \bra y {x} \leq f(z) - \bra yz\ \ \ \forall \, z\in X  \cr
\Iff &f(z) - \bra y {z} \ \ {\rm has\ a\ minimum\ point\ at\ \ } z=x.
}
$$
Define $-g(y)$ to be the minimum value, so that:
$$
f(x) + g(y) \ =\ \bra xy\qquad \forall\, (x,y)\in\partial f.
\eqno{(\HH.6)}
$$
Now
$$
\eqalign
{
x\in\partial g(y)  \Iff &g(y) + \bra x {w-y} \leq g(w) \ \ \ \forall \, w\in Y \cr
\Iff &g(y) -  \bra xy \leq g(w) - \bra xw\ \ \ \forall \,  w\in Y \cr
\Iff &g(w) - \bra xw\ \ {\rm has\ a\ minimum\ point\ at\ \ } w=y.
}
$$
The minimum value is $g(y)-\bra xy$, which equals $f(x)$.  This is the proof of the classical
fact that the Legendre transform is an involution.

Summarizing, if $f$ and $g$ correspond under the Legendre transform, then
$$
y\in \partial f \Iff x\in \partial g(y).
\eqno{(\HH.7)}
$$
That is, the multi-valued maps $F=\partial f$ and $G=\partial g$ are inverses of each other.

Note that $g$ convex, since $g(y)$ is the supremum of the 
family $\bra yx -f(x)$ of affine functions of $y$.

The Legendre transform $g$ of $f$ enjoys nicer properties
 than the convex function $f(x) \equiv r u(x) +\half |x|^2$.
We state more than required.

\Lemma {\HH.3}  {\sl
The Legendre transform $g$ of $f$ is a convex $C^1$ function on $X_\d$ with derivative
$Dg=G$.  If $G$ is differentiable at $y$ with first derivative $D_yG=B$, then $g$ 
is twice differentiable at $y$ with second derivative $D_y^2g =B$.
}
\pf
Since $\partial g = G$ and $G$ is single valued, we have $g\in C^1$ by Corollary  \BBB.3.
This enables us to apply the standard Mean Value Theorem.  Now assume that 
$G$ is differentiable at $y_0$ with $D_{y_0}G=B$.  We can assume $y_0=0$, 
$g(y_0)=0$,  and $G(y_0) \equiv D_{y_0}g=0$,  so that
$$
G(y)  -By \ =\ o(|y|).
\eqno{(\HH.8)}
$$
By the Mean Value Theorem applied to the function $\phi(y)\equiv
g(y)  - \half\bra {By}{y}$, 
 there exists $\x\in [0,y]$ such that 
$$
g(y)  - \half\bra {By}{y}  \ =\  \phi(y) \ =\  \bra { D_\x \phi} y\ =\ \bra{D_\x g  - B\x}{y}.
\eqno{(\HH.9)}
$$
By (\HH.8)
$$
D_\x g  -B\x \ =\ G(\x) -B(\x) \ =\ o(|\x|) \ =\ o(|y|).
$$
(This last equality is because $|\x|\leq |y|$.)  Therefore, the right hand side of (\HH.9) is $o(|y|^2)$.
This proves $\phi(y) = o(|y|^2)$ and hence $g$ is twice differentiable at 0 with $D^2_0g=B$.\qed 

\medskip

Now we are ready to prove the analogous lemma for $f$.

\Lemma{\HH.4} {\sl
Suppose that $G$ is differentiable at $y_0\in X_\d$ and 
let $B\equiv D_{y_0} G$ denote the derivative.  Assume 
that $x_0=G(y_0)$ is not a critical value of $G$.
Further assume that the convex function $f$ is differentiable at $x_0$, and hence $D_{x_0}f=y_0$.
Then the  function $f$ is twice differentiable at $x_0$ with second derivative $D_{x_0}^2f= B^{-1}$.}

\pf
We can assume that $x_0=0$ and that $f(0)= D_0f = 0$, by modifying $f$ by an affine function.
Since $f$ is differentiable at $0$, the subdifferential  $\partial f (0) = \{y_0\}  =\{D_{0}f\}=\{0\}$  by 
Lemma \BBB.2.
Since $y_0=0$ is not a critical point of $G$, the derivative $D_{0} G \equiv B$ is invertible.
Let $A\equiv B^{-1}$.
We must show that:
$$
f(x) -\half \bra{Ax}x \ =\ o \left(|x|^2\right).
\eqno{(\HH.10)}
$$

For $(x,y)\in\partial f$ the identity $f(x) + g(y) \ =\ \bra xy$ can be written as
$$
f(x) -\half \bra{Ax}x \ =\  \half \bra{By}y -g(y) +\half \bra {y-Ax} x + \half \bra{x-By} y. 
\eqno{(\HH.11)}
$$
We have by  (\HH.8) that
$$
x-By \ =\ G(y)-By \ =\ o(|y|).
\eqno{(\HH.8)'}
$$
Since $D_{y_0}^2g = B$, we have $g(y)-\half\bra{By}y = o(|y|^2)$. The remaining term of the 
RHS of (\HH.11) is also $o(|y|^2)$ since $|x|\leq |y|$ and 
$$
y-Ax \ =\ -A(x-By) \ =\ o(\|A\| |y|) \ =\ o(|y|).
\eqno{(\HH.12)}
$$
This proves that, with $x=G(y)$, 
$$
 f(x) -\half \bra{Ax}x \ =\  o\left(|y|^2\right).
\eqno{(\HH.13)}
$$
Finally for $(x,y)\in \partial f$, $|y| = |ABy| \leq \|A\|(|x-By|+|x|)$. 
 By (\HH.8)$'$ we have $|x-By| = o(|y|)$, so this proves $|y| = O(|x|)$.\qed
 
 \Remark{} Note that $|y|=O(|x|)$ combined with (\HH.12) yields $y-Ax = o(|x|)$.
 This is the statement that
 $$
 \lim_{\matrix{x\to x_0\cr y\in\partial f(x)}} {y-y_0-A(x-x_0)  \over |x-x_0|}\ =\ 0,
 \eqno{(\HH.14)}
$$
which says that the multi-valued function $\partial f$ is differentiable at $x_0$ with derivative $A$.

\medskip

Let $D$ denote the set of points in $G(X_\d)$ where the convex function $f$ is differentiable.
Let $N$ denote the set of points in $X_\d$ where $G$ is not differentiable.
Let $C$ denote the set of critical points for $G$ in $X_\d$ (where $G$ is differentiable
but the derivative of $G$ is singular).  Lemma \HH.4 applies to each point $x_0 = G(y_0)\in D$
with $y_0\notin N\cup C$.

\Lemma{\HH.5} {\sl
The set $D-G(N\cup C)$ has full measure in $G(X_\d)$.
}
\pf
As noted in (\HH.1) $D$ has full measure.
By Rademacher's Theorem  $N$ has measure zero.  Since $G$ is contractive,
$G(N)$ also has measure zero.  Finally, by the Lipschitz version of Sard's Theorem $G(C)$ has measure zero.\qed
\medskip

This proves that $f$ is twice differentiable on a set of full measure in $G(X_\d)$. 
Since $X_{2\d} \ss G(X_\d)$, the proof of Alexandrov's Theorem is complete.\qed

\vfill\eject
%\vskip .3in

%%%%%%%%%%%%%%%%%%%%%%%%%%%%%%%%%%%%%%%%%%%%%%%%%%%%
%%%%%%%%%%%%%%%%%%%%%%%%%%%%%%%%%%%%%%%%%%%%%%%%%%%%
%%%%%%%%%%%%%%%%%%%%%%%%%%%%%%%%%%%%%%%%%%%%%%%%%%%%
%%%%%%%%%%%%%%%%%%%%%%%%%%%%%%%%%%%%%%%%%%%%%%%%%%%%
%%%%%%%%%%%%%%%%%%%%%%%%%%%%%%%%%%%%%%%%%%%%%%%%%%%%

\centerline{\headfont Appendix  A.  A Quasi-Convexity Characterization of C$^{1,1}$.}
 \medskip 
 It is interesting that the condition that a function be $C^{1,1}$ is directly related to
 quasi-convexity, in fact it is equivalent to the function being simultaneously quasi-convex and quasi-concave.
 This was probably first observed by Hiriart-Urruty and Plazanet in [HP]. An alternate proof appeared in [E].
 For the benefit of the reader we include a proof here.

\Theorem {A.1}
$$
u \ \ {\rm is}  \ \l-C^{1,1} 
\qquad\iff\qquad
{\rm both}\ \ 
\pm u \ \ {\rm are}\ \ \l-{\rm quasi-convex}
$$

\pf 
($\Rightarrow$)  Suppose that $u$ is $\l$-$C^{1-1}$, i.e., $u\in C^1$ and $|D_xu -D_yu|\leq\l |x-y|$
for all $x,y$.  Set $f\equiv u+{\l\over 2}|x|^2$. Then 
$$
D_xf -D_yf\ =\ \l(x-y) + D_xu -D_yu,
$$
and hence
$$
\eqalign
{
\bra{D_xf -D_yf}{x-y} \ &=\ \l|x-y|^2 + \bra{D_xu -D_yu}{x-y}   \cr
&\geq \ \l|x-y|^2  - |D_xu -D_yu||x-y|   \cr
&= \ (\l|x-y|  - |D_xu -D_yu|)  |x-y|   \ \geq\ 0. 
}
$$
This form of monotonicity of $Df$ is one of the standard definitions of $f$ being convex.
The same proof works for $-u$

($\Leftarrow$ ) We state the converse as a proposition.

\Prop{A.2}
{\sl
If $u$ and $-u$ are $\l$-quasi-convex, then $u\in C^{1}$ and}
$$
|D_xu -D_yu|\ \leq \ \l|x-y|, \ \ \ {\rm i.e., \ } u\ \ {\rm is\ \ } \l- C^{1,1}.
$$
\pf
We first show that this is true if $u\in C^\infty$. 
Note that  $\pm u$ are $\l$-quasi-convex
\ \ $\iff$\ \ $D^2_x u +\l I\geq0$ and  $-D^2_x u +\l I\geq0$  for all $x$
\ \ $\iff$\ \ $-\l I \leq D^2_x u \leq \l I$ for all $x$.  By the Mean Value Theorem,
$D_xu -D_yu = (D^2_\x u)(x-y)$ for some $\x\in [x,y]$, and hence
$|D_xu -D_yu| \leq \l|x-y|$.  

In general, since the graph of $u+{\l\over 2}|x|^2$ has a supporting hyperplane from below and
 the graph of  $u-{\l\over 2}|x|^2$ has a supporting hyperplane from above, at every point, 
 the function $u$ is differentiable everywhere.  By partial continuity of the first derivative
 for quasi-convex functions (Lemma 1.3), we have $u\in C^1$.
 
 Now standard convolution $u^\e \equiv u * \vf_\e$ works just fine to complete the proof
 since $\pm u^\e$ is $\l$-quasi-convex by the next lemma, and the fact that
  $u\in C^1 \ \Rightarrow \ Du^\e \to Du$ locally uniformly.

\Lemma{A.3}
{\sl
$u$ is $\l$-quasi-convex \ \ $\Rightarrow$\ \ $u^\e \equiv u * \vf_\e$ is  $\l$-quasi-convex.
}
\pf
Suppose $u$ is  $\l$-quasi-convex, i.e., $f\equiv u+{\l\over 2}|x|^2$ is convex.
Standard convolution of $f$ with an approximate identity $\vf_\e$ based on $\vf$
(i.e., $\vf_\e(x) \equiv {1\over \e^n}\vf({x\over \e})$) yields $f^\e \equiv f * \vf_\e$ smooth and convex.
Note that $(|x|^2* \vf_\e) = |x|^2 +\bra ax+c$ preserves $|x|^2$ modulo an affine function,
since $\int|x+\e y|^2 \vf(y)\,dy = |x|^2+\e\bra ax+C\e^2$
where $\bra ax = 2\int \bra xy \vf(y) \, dy $ and
$C =\int|y|^2\vf(y)\,dy$.
Therefore, $D^2 f^\e = D^2 u^\e +\l I$, proving that each $u^\e$ is $\l$-quasi-convex.\qed

%\vfill\eject
\vskip .3in

%%%%%%%%%%%%%%%%%%%%%%%%%%%%%%%%%%%%%%%%%%%%%%%%%%%%
%%%%%%%%%%%%%%%%%%%%%%%%%%%%%%%%%%%%%%%%%%%%%%%%%%%%
%%%%%%%%%%%%%%%%%%%%%%%%%%%%%%%%%%%%%%%%%%%%%%%%%%%%
%%%%%%%%%%%%%%%%%%%%%%%%%%%%%%%%%%%%%%%%%%%%%%%%%%%%
%%%%%%%%%%%%%%%%%%%%%%%%%%%%%%%%%%%%%%%%%%%%%%%%%%%%

\centerline{\headfont   References.}
 \medskip

\item{[AA]}  G.   Alberti, and L. Ambrosio,
{\sl A geometrical approach to monotone functions in $\rn$}, 
 Math. Z.  {\bf 230} (1999), no. 2, 259-316.

\smallskip

\item{[Al]}  A. D. Alexandrov, {\sl Almost everywhere existence of the 
second differential of a convex function and properties of 
convex surfaces connected with it (in Russian)}, 
Lenningrad State Univ. Ann. Math.  {\bf 37}   (1939),   3-35.

\smallskip

\item{[CIL]}   M. G. Crandall, H. Ishii and P. L. Lions {\sl
User's guide to viscosity solutions of second order partial differential equations},  
Bull. Amer. Math. Soc. (N. S.) {\bf 27} (1992), 1-67.

 \smallskip

\item{[E]}  A. Eberhard, {\sl Prox-regularity and subjets}, in: A. Rubinov (Ed.), ÔOptimization and Related TopicsÕ, Applied Optimization Volumes, Kluwer Academic Publishers, Dordrecht, 2001, pp. 237Ð313.

\smallskip

\item{[F]}  H. Federer,   Geometric Measure Theory, Springer-Verlag, Berlin-Heidelberg, 1969.

\smallskip

\item {[HL$_1$]}   F. R. Harvey and H. B. Lawson, Jr,  {\sl  Dirichlet duality and the non-linear Dirichlet problem},    Comm. on Pure and Applied Math. {\bf 62} (2009), 396-443. ArXiv:math.0710.3991

\smallskip

\item {[HL$_2$]}  F. R. Harvey and H. B. Lawson, Jr.,   {\sl  The AE-Theorem and addition theorems
for quasi-convex functions}.   ArXiv:1309:1770.
\smallskip

\item{[HP]} J.-B. Hiriart-Urruty, Ph. Plazanet, {\sl MoreauÕs Theorem revisited}, Analyse Nonlin\'eaire (Perpignan, 1987), Ann. Inst. H. Poincar\'e {\bf 6} (1989), 325Ð338.

 \smallskip

\item {[S]}  Z. Slodkowski, {\sl  The Bremermann-Dirichlet problem for $q$-plurisubharmonic functions},
Ann. Scuola Norm. Sup. Pisa Cl. Sci. (4)  {\bf 11}    (1984),  303-326.

\smallskip

\item{[J]}    R. Jensen,    {\sl  The maximum principle for viscosity solutions of fully nonlinear
second order partial differential equations},   Arch. Ratl. Mech. Anal. {\bf 101}  (1988),   1-27.

\smallskip

%\item{[M]}    F. Mignot,    {\sl  Contr\^ol optimal dans les in\'equations variationelles elliptiques},   %Ratl. J. Funct. Anal. {\bf 22}  (1976),   130-185.

%\smallskip

\end